\documentclass[preprint,12pt]{elsarticle}
\usepackage{amsfonts}
\usepackage{amsmath}
\usepackage{mathrsfs}
\usepackage[applemac]{inputenc}
\usepackage{amssymb}
\usepackage{times, ulem, amsfonts}
\usepackage{mathpazo, amsmath, amssymb, amsthm,color}
\makeatletter
\def\ps@pprintTitle{%
   \let\@oddhead\@empty
   \let\@evenhead\@empty
   \let\@oddfoot\@empty
   \let\@evenfoot\@oddfoot
}
\makeatother

\def\q{{\mathbf Q}}
\def\v{{\mathcal V}}
\def\d{{\mathcal D}}
\def\u{\mathrm u}
\def\p{{\mathbf P}}
\def\I{{\mathbf I}}
\newcommand{\R}{{\mathbb{R}}}
\newcommand{\Z}{{\mathbb{Z}}}
\newcommand\diverg{\mathop{\mbox{\rm div}}}
\def\Rn{({\mathbb R}^n)}
\newtheorem{theorem}{Theorem}[section]

\newtheorem{remark}[theorem]{Remark}
\numberwithin{equation}{section}

\begin{document}
\title{ Global large solutions and incompressible limit for the compressible flow of liquid crystals}

  \author{Xiaoping Zhai}
\ead{zhaixp@szu.edu.cn}

 \author{Zhi-Min Chen\corref{cor1}}
 \ead{zmchen@szu.edu.cn}
 \address{School  of Mathematics and Statistics, Shenzhen University, Shenzhen 518060, China}
\cortext[cor1]{Corresponding author}

\baselineskip=24pt

\begin{abstract}
The present paper is dedicated to the existence of global large solutions  and incompressible limit for the compressible flow of liquid crystals under the assumption on almost constant density and large bulk viscosity. The result is based on the Fourier analysis and involves so-called critical homogeneous Besov spaces.
\end{abstract}

\maketitle


\section{Introduction}
Consider  the compressible flow of liquid crystals in $\R^n (n=2, 3)$:
\begin{eqnarray}\label{m}
\left\{\begin{aligned}
&\partial_t \rho \!+\! \diverg(\rho v) = 0\, ,\\
 &     \partial_t ( \rho v ) \!+\! \diverg ( \rho v \!\otimes\! v ) \!-\!\mu\Delta v\!-\!(\mu\!+\!\lambda)\nabla\diverg v\!+\! \nabla P( \rho) = \!-\!\diverg ( \nabla d\!\odot\! \nabla d\!-\!\tfrac12|\nabla d|^2\I)\, ,\\
   &   \partial_td \!+\!v\cdot\nabla d= \Delta d  \!+\! |\nabla d|^2d\, ,\\
    &(\rho,v,d)|_{t=0}=(\rho_0,v_0,d_0),
\end{aligned}\right.
\end{eqnarray}
where
$\rho$  is the  density field, $v$ is the  velocity field,  $P$ is  the pressure field and $d$ valued in $\mathbb{S}^{n-1}$ represents
the macroscopic average of the nematic liquid crystal orientation field. The symbol $\nabla d \odot\nabla d$, which exhibits the property of the anisotropy of
the material, denotes the matrix $(\partial_id\cdot \partial_jd)_{n\times n}$.
The parameters  $\mu$ and $\lambda$   are shear viscosity and bulk viscosity coefficients respectively satisfying the physical condition
 $$ \mu >0,\,\,\, \nu \triangleq \lambda+2\mu > 0$$
  which  ensures ellipticity property
for the Lam\'{e} operator $\mu\Delta + (\lambda + \mu)\nabla\diverg$ .

 In the present paper, we are interested in  \eqref{m} under the  isentropic flow assumption:
\begin{equation*}
P(\rho)= A\rho^\gamma\quad{\rm with}\quad
A>0,\ \gamma>1. 
\end{equation*}
The flow governed  by  (\ref{m})  is basically a coupling of a compressible Navier-Stokes flow  and a parabolic heat flow.
When $d\equiv {\it constant}$, Eq. (\ref{m}) reduces to the well-known
Navier-Stokes equations for the compressible isentropic flow,  which has been studied by many researchers (see \cite{danchin2007cpde,danchin2000,danchin2018,Hoff,huangxiangdi2012} and the references therein).

The existence of global weak solutions
of  \eqref{m} in $\R^2$  was obtained by  Jiang {\it et al.} \cite{jjw}  under the assumption that the image of the initial data of $d$ is contained in $\mathbb{S}^2_+$. For the compressible nematic liquid crystal flow  with $d$ having  variable degree of orientations, the global existence of weak solutions in $\R^3$  has been obtained by \cite{liu2013} and \cite{WangA2012}. Recently, inspired by the work of \cite{linfanghua2016} on  incompressible flow, the corresponding global finite energy weak solutions to \eqref{m} was proved in \cite{linfanghua2015}.
The local existence of strong solutions in $\R^3$ has been  studied by \cite{huang2012} and \cite{huang2012+1}.
There are also many results 
 dealing  with the incompressible limit of compressible nematic liquid crystal flow.
Hao and Liu \cite{hao} investigated  the  incompressible limit  in the whole space and a bounded domain  under the  Dirichlet boundary condition.
Ding {\it et al.} \cite{ding1} studied the incompressible limit in bounded domain   with periodic boundary condition.
Wang and Yu \cite{wang2014} proved the incompressible limit for weak solutions in a bounded domain. For more study about the incompressible limit problem, one can refer to \cite{bai1,bai2,lions,ouraobing} and the references therein.

 Recently, Danchin and Mucha \cite{danchin2018}
established a technical method for the  global solutions of incompressible equations as the limit ($\lambda \to \infty$) of those of  the compressible isentropic Navier-Stokes
equations supplemented with arbitrary large initial velocity  and almost constant density.
Motivated by \cite{danchin2018},  we  study the global  large solutions and incompressible limit for  \eqref{m}
with the stronger nonlinear external forcing terms   $\diverg ( \nabla d\odot \nabla d-\tfrac12|\nabla d|^2\I)$ and $|\nabla d|^2d$.
 With the development of the method in \cite{danchin2018},
 we obtain global large strong solutions to \eqref{m} with respect to the bulk viscosity coefficient   $\lambda$ sufficiently large, any shear viscosity coefficient   $\mu>0$  and the initial density $\rho_0$ being  very close to a positive constant.

The limit of  \eqref{m}  as $\lambda\to +\infty$ is to be shown to  the following system:
\begin{align}\label{mm}
\left\{\begin{aligned}
&\partial_{t} \v  -\mu\Delta \v  +\v  \cdot \nabla \v  +\nabla \Pi=-\diverg ( \nabla \d\odot \nabla  \d),\\
&\partial_{t}  \d+\v  \cdot \nabla   \d=\Delta \d+|\nabla  \d|^2 \d,\\
&\mathrm{div}\v  =0,\\
&\v  |_{t=0}=\p v_0,\quad \d|_{t=0}=d_0,
\end{aligned}\right.
\end{align}
with the projector $\p=\I-\q\triangleq \I-\nabla\Delta^{-1}\diverg$.

Now, we state the main result of the present paper:
\begin{theorem}\label{dingli}
 Let  $\nu >0$ and $(v_0,d_0)  \in \dot B^{-1+\frac n2}_{2,1}(\R^n)\times  \dot B^{\frac n2}_{2,1}(\R^n)$ for $n=2,3$. Let  $$a_0=\rho_0-1\in\dot B^{-1+\frac n2}_{2,1}\Rn\cap \dot B^{\frac n2}_{2,1}(\R^n).$$
 Suppose that \eqref{mm} with the initial data $\v  _0=\p v_0$ and $\d|_{t=0}=d_0$
 generates a unique global solution $(\v , \d) \in C(\R^+;\dot B^{-1+\frac n2}_{2,1}\Rn)\times C(\R^+;\dot B^{\frac n2}_{2,1}\Rn)$ such that
\begin{align}\label{f3a}
&  \|\v  \|_{L^\infty(0,T;\dot B^{-1+\frac{n}{2}}_{2,1})}+ \|\d \|_{L^\infty(0,T;\dot B^{\frac{n}{2}}_{2,1})}\nonumber\\
&\quad+\|\v  _t\|_{L^1(0,T;\dot B^{-1+\frac{n}{2}}_{2,1})}+\| \v  \|_{L^1(0,T;\dot B^{1+\frac{n}{2}}_{2,1})}+\| \d \|_{L^1(0,T;\dot B^{2+\frac{n}{2}}_{2,1})}
\leq M
\end{align}
with respect to  a constant bound  $M$. For $\mu>0$ and a constant $C>0$, assume that the parameter  $\nu$ is subject to the inequality
\begin{align}\label{new}
 C \exp(CM+CM^4)\left(\|a_0\|_{\dot B^{-1+\frac n2}_{2,1}}+\nu\|a_0\|_{\dot B^{\frac n2}_{2,1}}+\|\q v_0\|_{\dot B^{-1+\frac n2}_{2,1}}+M^2
 +\mu^2 \right)\leq \sqrt{\mu\nu}.
 \end{align}
 Then   \eqref{m} has a unique global  solution $(\rho,v,d)$  such that
\begin{align*}
&v \in C(\R^+;\dot B^{-1+\frac n2}_{2,1}\Rn),\\
& v_t,\nabla^2v\in L^1(0,\infty; \dot B^{-1+\frac n2}_{2,1}\Rn),   \nonumber\\
&d \in C(\R^+;\dot B^{\frac n2}_{2,1}\Rn)\cap L^1(0,\infty; \dot B^{2+\frac n2}_{2,1}\Rn),\nonumber\\ 
&a\triangleq\rho -1 \in C(\R^+; \dot B^{-1+\frac n2}_{2,1}\Rn\cap \dot B^{\frac n2}_{2,1}\Rn) \cap
  L^2(0,\infty; \dot B^{\frac n2}_{2,1}\Rn),
\end{align*}
and
\begin{align}\label{}
 &\|\q v\|_{L^\infty(0,\infty; \dot B^{-1+\frac n2}_{2,1})}+
  \| \nu\nabla^2\q v\|_{L^1(0,\infty; \dot B^{-1+\frac n2}_{2,1})}
    + \|a\|_{L^\infty(0,\infty; \dot B^{-1+\frac n2}_{2,1})} +
      \nu\|a\|_{L^\infty(0,\infty; \dot B^{\frac n2}_{2,1})}\nonumber\\
    &\quad\leq   C\exp(CM+CM^4)\bigl(\|a_0\|_{\dot B^{-1+\frac n2}_{2,1}} +\nu \|a_0\|_{\dot B^{\frac n2}_{2,1}}+\|\q v_0\|_{\dot B^{-1+\frac n2}_{2,1}} +M^2+\mu^2\bigr).
        \end{align}
In addition,  if $\rho_0=1$, then  $(\rho,v,d)\to(1,\v  ,\d )$  as $\nu\to \infty$ in the following sense:
 \begin{align}\label{}
 &\sqrt{\mu\nu^{-1}}\big\|\rho-1\big\|_{L^\infty(0,\infty; \dot B^{\frac n2}_{2,1})}
+  \big\|\p v -\v  \big\|_{L^\infty(0,\infty; \dot B^{-1+\frac n2}_{2,1})}  +  \big\|d -\d \big\|_{L^\infty(0,\infty; \dot B^{\frac n2}_{2,1})}   \nonumber\\
&\quad\quad+\big\|(\p v_t -\v  _t,\mu\nabla^2 (\p v -\v  ))\big\|_{L^1(0,\infty; \dot B^{-1+\frac n2}_{2,1})}  +  \big\| d-\d \big\|_{L^1(0,\infty; \dot B^{2+\frac n2}_{2,1})}\nonumber\\
&\quad   \leq C\sqrt{\mu\nu^{-1}}.
  \end{align}
\end{theorem}
 \begin{remark}
 Without
loss of generality, we will fix   $\mu = 1$ and so $\nu=\lambda +2$ in  the proof of the result.
\end{remark}
The rest of the paper unfolds as follows. In Section 2,  we recall some  basic facts of Littlewood-Paley
theory and some basic estimates.
The proof of the main result is detailed in  Section 3, where the analysis is separated into three subsections for  incompressible part estimates,  compressible part estimates  and   combination estimates.
A blow-up  criterion  is examined in the Appendix.



\section{Preliminaries } Denote  by $C$  a generic constant, which may vary from line to line.
 The symbol $\mathcal{F}$ represents  the Fourier transform in $\mathcal{S}'(\R^n)$, the space of tempered distributions in $\R^n$. Let $\varphi$ be a nonnegative smooth function supported in an annulus  of $\R^n$ so that
\begin{align*}
 \sum_{j \in \mathbb{Z}} \varphi_j(\cdot)=1 \ \mbox{ in } \ \R^n\setminus \{0\} \ \mbox{ for } \,\,\varphi_j(\cdot)= \varphi(2^{-j}\cdot).
 \end{align*}
Therefore the homogeneous  dyadic blocks are defined as
\begin{align*}
&\dot{\Delta}_ju=\mathcal{F}^{-1}( \varphi_j\mathcal{F}u).
\end{align*}
Denote by $\mathcal{S}_h^{'}({\mathbb R} ^n)$ the space of $u\in \mathcal{S}'(\R^n)$ such that
 the decomposition
$$
u=\sum_{ j\in \mathbb{Z} }\dot{ \Delta }_{j}u,\ \ \
$$
holds true. Hence we have the
homogeneous Besov space
 $$\dot{B}_{2,1}^s(\mathbb{R}^n)=\left\{u\in \mathscr{S}'_h(\mathbb{R}^n)\left|\, \|u\|_{\dot{B}_{2,1}^s}\triangleq \sum_{j\in\Z}2^{js}\|\dot{\Delta}_ju\|_{L^2}<\infty\right.\right\} \mbox{ for } s \in \mathbb{R}.$$

For convenience, we use the symbols
\begin{align*}
&L^p_{T}(\dot{B}_{2,1}^s(\mathbb{R}^n))\triangleq L^p(0,T; \dot{B}_{2,1}^s(\mathbb{R}^n)),\\
&{\widetilde{L}^p_{T}(\dot{B}_{2,1}^s(\mathbb{R}^n))}\triangleq \left\{u: [0,T ]\mapsto \mathcal{S}'(\mathbb{R}^n)\left|\,
\|u\|_{\widetilde{L}^p_{T}(\dot{B}_{2,1}^s)}\triangleq \sum_{j\in \mathbb{Z}}2^{js}
\|\dot{\Delta}_ju\|_{L^p(0,T;L^2(\R^n))}<\infty\right.\right\}.
\end{align*}

It is readily seen  that  there holds the  interpolation inequality:
\begin{align}\label{inter}
\|u\|_{\widetilde{L}^p_{T}(\dot{B}_{2,1}^s)}\lesssim
\|u\|^\theta_{\widetilde{L}^{p_1}_{T}(\dot{B}_{2,1}^{s_1})}
\|u\|^{1-\theta}_{\widetilde{L}^{p_2}_{T}(\dot{B}_{2,1}^{s_2})}, \, \frac{1}{p}=\frac{\theta}{p_1}+\frac{1-\theta}{p_2}, \, s=\theta s_1+(1-\theta)s_2
\end{align}
for $0<s_1<s_2,$ $\theta\in[0,1]$ and $p_1,\, p_2\in[1,+\infty]$.

We will also repeatedly use the  following Bernstein inequality:
\begin{align}
&& C^{-1}\sigma ^k\|u\|_{L^q}\le\|\nabla^k u\|_{L^q}
\le C\sigma ^{k+\frac np-\frac nq}\|u\|_{L^p} \mbox{ when } \mathrm{Supp} \,\mathcal{F}u\subset\sigma  \mathcal{C}\label{B2222}
\end{align}
for $1\le p \le q\le \infty$, $k\in \Z$,   $\mathcal{C}$ an annulus  of $\mathbb{R}^n$ and $C$ a generic constant independent of the scale parameter $\sigma >0$.

Moreover, we will use the following pointwise product law \cite[Lemma 2.7]{zhaixiaoping}
\begin{align}\label{daishu}
\|uv\|_{\dot{B}_{2,1}^{s_1+s_2 -\frac{n}{2}}}\lesssim \|u\|_{\dot{B}_{2,1}^{s_1}}\|v\|_{\dot{B}_{2,1}^{s_2}},\,\,\, s_1\leq \frac{n}{2}, \,\, \,s_2\leq \frac n2,\,\,\,s_1+s_2>0
\end{align}
and the commutator estimate \cite[Remark 2.102]{bcd}
\begin{align}\label{jiaohuanzi}
\sum_{j\in \mathbb{Z}}2^{js}\left\|[u\cdot \nabla ,\dot{\Delta}_j]v\right\|_{L^2}\lesssim  \|\nabla u\|_{\dot{B}_{2,1}^{\frac{n}{2}}}\|v\|_{\dot{B}_{2,1}^{s}}, \  s= \frac n2, \ -1+\frac n2
\end{align}
for the commutator $[u\cdot \nabla ,\dot{\Delta}_j]= u\cdot \nabla \dot\Delta_j - \dot\Delta_j (u\cdot \nabla )$.

\section{Proof of Theorem \ref{dingli}}
In this section,
 the global wellposedness of \eqref{m} will be established by extending
a local solution through  global in time a priori estimates and a continuous criterion. In the rest of this paper, we focus on the global in time a priori estimates.

For the local solution initially from the  data $a_0, v_0, d_0$ given  in Theorem \ref{dingli},
 we can simply follow  a local solution  method from Danchin  \cite{danchin2007cpde} or Hu and Wu \cite{huxianpeng} to obtain the  local wellposedness of system \eqref{m}:
\begin{theorem}\label{dingli2}
 Let
 $v_0 \in \dot B^{-1+\frac n2}_{2,1}(\R^n)$,   $d_0 \in \dot B^{\frac n2}_{2,1}(\R^n)$  and $a_0=\rho_0-1 \in\dot B^{\frac n2}_{2,1}(\R^n)$ with $n=2,\,3$ and
$1+\inf_{x\in\R^n}a_0(x)>0$. Then there exists a positive time $T^\ast$ so that \eqref{m} has a unique solution $(a,v,d)$ on $[0, T^\ast )$ so that for any $T<T^\ast $,
\begin{align*}
&a\in C([0,T];\dot{B}_{2,1}^{{\frac n2}}\Rn)\cap
\widetilde{L}^\infty_T(\dot{B}_{2,1}^{{\frac n2}}\Rn),\\
&v\in C([0,T];\dot{B}_{2,1}^{{\frac n2-1}}\Rn)\cap
\widetilde{L}^\infty_T(\dot{B}_{2,1}^{{\frac n2-1}}\Rn)\cap
L^1_T(\dot{B}_{2,1}^{{\frac n2+1}}\Rn),\\
&d\in C([0,T];\dot{B}_{2,1}^{{\frac n2}}\Rn)\cap
\widetilde{L}^\infty_T(\dot{B}_{2,1}^{{\frac n2}}\Rn)\cap
L^1_T(\dot{B}_{2,1}^{{\frac n2+2}}\Rn).
\end{align*}
Moreover,
if $T^\ast$ is the maximal time of existence to this solution and $T^\ast<\infty$, then one has
\begin{align}\label{baopobiaozhun}
\int_0^{T^\ast}(\|\nabla v\|_{\dot{B}_{2,1}^{\frac{n}{2}}}+\|d\|^2_{\dot{B}_{2,1}^{1+\frac{n}{2}}})dt=\infty.
\end{align}
\end{theorem}
Thus the proof of  the  existence and uniqueness assertion of Theorem \ref{baopobiaozhun} is omitted. However,   the demonstration of the blow-up criterion \eqref{baopobiaozhun} is shown in Appendix.

We also assume that the system \eqref{mm} generates a global large solution with initial data $\p v_0\in\dot{B}_{2,1}^{-1+\frac{n}{2}}(\R^n) $ and  $d_0\in\dot{B}_{2,1}^{\frac{n}{2}}(\R^n) $. In fact, this assumption can be verified from \cite{liuqiao} and \cite{zhaixiaoping}. More precisely, if the initial data  $\p v_0$ and $d_0$ satisfy a smallness condition, we can get that   \eqref{mm} has a unique global solution $(\v,\d)$ such that
\begin{align*}\label{}
&\v\in C([0,T];\dot{B}_{2,1}^{{\frac n2-1}}\Rn)\cap
L^1_T(\dot{B}_{2,1}^{{\frac n2+1}}\Rn),\nonumber\\
&\v_t,\nabla^2\v\in
L^1_T(\dot{B}_{2,1}^{{-1+\frac n2}}\Rn),\nonumber\\
&\d\in C([0,T];\dot{B}_{2,1}^{{\frac n2}}\Rn)\cap
L^1_T(\dot{B}_{2,1}^{{\frac n2+2}}\Rn).
\end{align*}
Moreover, we  assume
\begin{align*}
{V}_n(T)\triangleq  \|\v  \|_{L^\infty_T(\dot B^{-1+\frac{n}{2}}_{2,1})}\!+\! \|\d \|_{L^\infty_T(\dot B^{\frac{n}{2}}_{2,1})}\!+\!\|\v  _t\|_{L^1_T(\dot B^{-1+\frac{n}{2}}_{2,1})}
\!+\!\| \v  \|_{L^1_T(\dot B^{1+\frac{n}{2}}_{2,1})}\!+\!\| \d \|_{L^1_T(\dot B^{2+\frac{n}{2}}_{2,1})}\le M.
\end{align*}

With the aid of Theorem \ref{dingli2}, it remains to show global in time a priori estimates.

The proof is divided into three subsections. By applying   the operators $\p$ and $\q$ to decompose  the system \eqref{m} into the incompressible part and the compressible part, the estimates for the incompressible part is given in the first subsection, while the   estimates of the compressible part is shown in the second subsection. Then the combination of the two parts estimates for producing desired global in time  energy estimates  is given in the last subsection.

\subsection{Estimates for the incompressible part }
Let
$$ \u = v-\v  \ \mbox{ and } \ \delta=d-\d .$$
Applying $\p$ to the second equation in \eqref{m} and the first equation in \eqref{mm} respectively, then making  the difference between the equations of $(\p v, d)$ and those of $(\v ,\d)$, we get the incompressible component
\begin{eqnarray}\label{i}
\left\{\begin{aligned}
&(\p \u )_t+\p((\u +\v  )\cdot\nabla \p \u )-\Delta\p \u =-\p H_1,\\
&\delta_t -\Delta \delta+(\u +\v  )\cdot\nabla\delta=G(\delta,\d )-\u \cdot\nabla \d ,\\
&\p \u |_{t=0}=0,\quad \delta|_{t=0}=0,
\end{aligned}\right.
\end{eqnarray}
with
\begin{align*}
G(\delta,\d )=&|\nabla \delta|^2\delta+|\nabla \d |^2\delta+|\nabla \delta|^2\d +2\nabla\delta\cdot\nabla\d \cdot\delta+2\nabla\delta\cdot\nabla\d \cdot\d ,\\
H_1=&\underbrace{a\left(\v  _t+\p \u _t+(\q \u _t+\nabla a)\right)}_{H_1^{(1)}}-\underbrace{(1+a)\p \u \cdot\nabla(\v  +\q \u )}_{H_1^{(2)}}
\\
&+\underbrace{(1+a)(\v  \cdot\nabla\q \u +\q \u \cdot\nabla \v  )}_{H_1^{(3)}}+\underbrace{a(\u +\v  )\cdot\nabla \p \u }_{H_1^{(4)}}
\nonumber\\
&+\underbrace{a (\q \u \cdot\nabla \q \u + \v  \cdot\nabla \v  )}_{H_1^{(5)}}+\underbrace{\diverg(\nabla \delta\odot \nabla \delta+\nabla\delta\odot \nabla \d +\nabla\d \odot \nabla \delta)}_{H_1^{(6)}}.
\end{align*}
Similarly, applying $\q$ to the second equation in \eqref{m} and assuming $P'(1)=1$, we  get the compressible component
\begin{eqnarray}\label{c}
\left\{\begin{aligned}
&a_t+\diverg(a\u )+\diverg\q \u  +\v  \cdot\nabla a =0,\\
&(\q \u )_t+\q((\u +\v  )\cdot\nabla \q \u )-\nu\Delta\q \u +\nabla a=-\q H_2,\\
&a|_{t=0}=a_0,\quad \q \u |_{t=0}=\q v_0,
\end{aligned}\right.
\end{eqnarray}
with
\begin{align*}
H_2=&\  a\left(\v  _t+\p \u _t+(\q \u _t+\nabla a)\right)+\diverg ( \nabla (\delta+\d )\odot \nabla (\delta+\d )-\tfrac12|\nabla (\delta+\d )|^2\I)\nonumber\\
&+(1+a)(\u +\v  )\cdot\nabla \p \u +(1+a)(\u +\v  )\cdot\nabla \v  +a(\u +\v  )\cdot\nabla \q \u  +(k(a)-a)\nabla a,
\end{align*}
and
$$ k(a)=P'(1+a)-P'(1)=P'(1+a)-1.$$

For  the estimates of  the incompressible part, we  use the  notation:
\begin{align}
X_n(T)\triangleq &\big\|(\q \u , a,\nu\nabla a)\big\|_{L^\infty_T(\dot{B}_{2,1}^{-1+\frac{n}{2}} )},\label{new1}\\
 Y_n(T)\triangleq &\big\|(\q \u _t+\nabla a,\nu\nabla^2\q \u , \nu\nabla^2a^\ell, \nabla a ^h)\big\|_{L^1_T(\dot{B}_{2,1}^{-1+\frac{n}{2}} )},\\
Z_n(T)\triangleq &\big\|\p \u \big\|_{L^\infty_T(\dot{B}_{2,1}^{-1+\frac{n}{2}} )}+\big\|\delta\big\|_{L^\infty_T(\dot{B}_{2,1}^{\frac{n}{2}} )}, \\
W_n(T)\triangleq &\big\|(\p \u _t,\nabla^2\p \u )\big\|_{L^1_T(\dot{B}_{2,1}^{-1+\frac{n}{2}} )}+\big\|\delta\big\|_{L^1_T(\dot{B}_{2,1}^{2+\frac{n}{2}} )}.
\end{align}


Applying $\dot{\Delta_j}$ to the second equation in \eqref{i} gives
\begin{align}\label{B1}
\partial_t\dot{\Delta_j}\delta -\Delta \dot{\Delta_j}\delta+(\u +\v  )\cdot\nabla\dot{\Delta_j}\delta=[(\u +\v  )\cdot\nabla,\dot{\Delta_j}]\delta+\dot{\Delta_j}G(\delta,\d )-\dot{\Delta_j}(\u \cdot\nabla \d ).
\end{align}
Taking $L^2$ inner product with $\dot{\Delta_j}\delta$ and using integrating by part imply that
\begin{align}\label{B2}
\frac12\frac{d}{dt}\|\dot{\Delta_j}\delta\|_{L^2}^2+2^{2j}\|\dot{\Delta_j}\delta\|_{L^2}^2\lesssim&\|\diverg(\u +\v  )\|_{L^\infty}\|\dot{\Delta_j}\delta\|_{L^2}^2+\big(\|[\dot{\Delta_j},(\u +\v  )\cdot\nabla]\delta\|_{L^2}\nonumber\\
&+\|\dot{\Delta_j}G(\delta,\d )\|_{L^2}+\|\dot{\Delta_j}(\u \cdot\nabla \d )\|_{L^2}\big)\|\dot{\Delta_j}\delta\|_{L^2}.
\end{align}
By (\ref{jiaohuanzi}), we have
\begin{align}\label{B3}
\sum_{j\in\Z}2^{\frac n2j}\|[\dot{\Delta_j},(\u +\v  )\cdot\nabla]\delta\|_{L^2}\lesssim\|(\nabla \u ,\nabla \v  )\|_{\dot{B}_{2,1}^{\frac{n}{2}}}\|\delta\|_{\dot{B}_{2,1}^{\frac{n}{2}}}.
\end{align}
Inserting \eqref{B3} into \eqref{B2} and employing the embedding property $\dot{B}_{2,1}^{\frac{n}{2}}(\R^n)\hookrightarrow L^\infty(\R^n)$, we have
\begin{align}\label{A1}
\|\delta\|_{L^\infty_T(\dot{B}_{2,1}^{\frac{n}{2}} )}+\|\delta\|_{L^1_T(\dot{B}_{2,1}^{2+\frac{n}{2}} )}
\lesssim&\int_0^T\|(\nabla \u ,\nabla \v  )\|_{\dot{B}_{2,1}^{\frac{n}{2}}}\|\delta\|_{\dot{B}_{2,1}^{\frac{n}{2}}}dt\nonumber\\
&+\int_0^T\|G(\delta,\d )\|_{\dot{B}_{2,1}^{\frac{n}{2}}}dt+\int_0^T\|\u \cdot\nabla \d \|_{\dot{B}_{2,1}^{\frac{n}{2}}}dt.
\end{align}
By \eqref{daishu}, we have
\begin{align}\label{B5}
\int_0^T\|\u \cdot\nabla \d \|_{\dot{B}_{2,1}^{\frac{n}{2}}}dt
&\le C\int_0^T(\|\p \u \|_{\dot{B}_{2,1}^{\frac{n}{2}}}+\|\q \u \|_{\dot{B}_{2,1}^{\frac{n}{2}}})\| \d \|_{\dot{B}_{2,1}^{1+\frac{n}{2}}}dt\nonumber\\
&\le C\int_0^T(\|\p \u \|_{\dot{B}_{2,1}^{-1+\frac{n}{2}}}^\frac12\|\p \u \|_{\dot{B}_{2,1}^{1+\frac{n}{2}}}^\frac12+\|\q \u \|_{\dot{B}_{2,1}^{-1+\frac{n}{2}}}^\frac12\|\nabla^2\q \u \|_{\dot{B}_{2,1}^{-1+\frac{n}{2}}})^\frac12)\| \d \|_{\dot{B}_{2,1}^{1+\frac{n}{2}}}dt\nonumber\\
&\le
\varepsilon\int_0^T\|\p \u \|_{\dot{B}_{2,1}^{1+\frac{n}{2}}}dt+C\int_0^T\| \d \|_{\dot{B}_{2,1}^{1+\frac{n}{2}}}^2\|\p \u \|_{\dot{B}_{2,1}^{-1+\frac{n}{2}}}dt\nonumber\\
&+\nu^{-\frac12}\left(\int_0^T\|\q \u \|_{\dot{B}_{2,1}^{-1+\frac{n}{2}}}\|\nu\nabla^2\q \u \|_{\dot{B}_{2,1}^{-1+\frac{n}{2}}}dt\right)^\frac12
\left(\int_0^T\| \d \|_{\dot{B}_{2,1}^{1+\frac{n}{2}}}^2dt\right)^\frac12.
\end{align}
Similarly, it follows from \eqref{daishu}, the H\"older inequality, interpolation inequality \eqref{inter} and the Young inequality that
\begin{align}\label{A2}
\int_0^T&\|G(\delta,\d )\|_{\dot{B}_{2,1}^{\frac{n}{2}}}dt\nonumber\\
\lesssim&\int_0^T\big\||\nabla \delta|^2\delta\!+\!|\nabla \d |^2\delta\!+\!2\nabla\delta\cdot\nabla\d \cdot\delta\big\|_{\dot{B}_{2,1}^{\frac{n}{2}}}dt 
\!+\!\int_0^T\big\||\nabla \delta|^2\d \!+\!2\nabla\delta\cdot\nabla\d \cdot\d \big\|_{\dot{B}_{2,1}^{\frac{n}{2}}}dt\nonumber\\
\lesssim&\int_0^T(\|\delta\|_{\dot{B}_{2,1}^{1\!+\!\frac{n}{2}}}^2\!+\!\|\d \|_{\dot{B}_{2,1}^{1\!+\!\frac{n}{2}}}^2)\|\delta\|_{\dot{B}_{2,1}^{\frac{n}{2}}}dt 
\!+\!\int_0^T\|\d \|_{\dot{B}_{2,1}^{\frac{n}{2}}}\|\delta\|_{\dot{B}_{2,1}^{1\!+\!\frac{n}{2}}}^2\!+\!\|\d \|_{\dot{B}_{2,1}^{\frac{n}{2}}}\|\d \|_{\dot{B}_{2,1}^{1\!+\!\frac{n}{2}}}\|\delta\|_{\dot{B}_{2,1}^{1\!+\!\frac{n}{2}}}dt\nonumber\\
\lesssim&\int_0^T(\|\delta\|_{\dot{B}_{2,1}^{1\!+\!\frac{n}{2}}}^2\!+\!\|\d \|_{\dot{B}_{2,1}^{1\!+\!\frac{n}{2}}}^2)\|\delta\|_{\dot{B}_{2,1}^{\frac{n}{2}}}dt\nonumber\\
&\quad\!+\!\int_0^T\|\d \|_{\dot{B}_{2,1}^{\frac{n}{2}}}\|\delta\|_{\dot{B}_{2,1}^{2\!+\!\frac{n}{2}}}\|\delta\|_{\dot{B}_{2,1}^{\frac{n}{2}}}\!+\!\|\d \|_{\dot{B}_{2,1}^{\frac{n}{2}}}\|\d \|_{\dot{B}_{2,1}^{1\!+\!\frac{n}{2}}}\|\delta\|_{\dot{B}_{2,1}^{\frac{n}{2}}}^\frac12 \|\delta\|_{\dot{B}_{2,1}^{2\!+\!\frac{n}{2}}}^\frac12dt\nonumber\\
\lesssim&\!\!\int_0^T\!\!\!\|\delta\|_{\dot{B}_{2,1}^{2\!+\!\frac{n}{2}}}\frac{dt}2\!+\!\!\!\int_0^T\!\!\!\!\Big(\|\delta\|_{\dot{B}_{2,1}^{1\!+\!\frac{n}{2}}}^2\!+\!\|\d \|_{\dot{B}_{2,1}^{1\!+\!\frac{n}{2}}}^2 
 \!+\!\|\d \|_{\dot{B}_{2,1}^{\frac{n}{2}}}\|\delta\|_{\dot{B}_{2,1}^{2\!+\!\frac{n}{2}}}\!+\!\|\d \|_{\dot{B}_{2,1}^{\frac{n}{2}}}^2 \|\d \|_{\dot{B}_{2,1}^{1\!+\!\frac{n}{2}}}^2\Big)\|\delta\|_{\dot{B}_{2,1}^{\frac{n}{2}}}dt.
\end{align}
Let
\begin{align*}
 H_3(\p \u ,\q \u ,\v  ,\delta,\d )=&\|\nabla \p \u \|_{\dot{B}_{2,1}^{\frac{n}{2}}}+\|\nabla \q \u \|_{\dot{B}_{2,1}^{\frac{n}{2}}}+\|\nabla \v  \|_{\dot{B}_{2,1}^{\frac{n}{2}}}+\|\delta\|_{\dot{B}_{2,1}^{1+\frac{n}{2}}}^2\\
 &+\|\d \|_{\dot{B}_{2,1}^{1+\frac{n}{2}}}^2  +\|\d \|_{\dot{B}_{2,1}^{\frac{n}{2}}}\|\delta\|_{\dot{B}_{2,1}^{2+\frac{n}{2}}}+\|\d \|_{\dot{B}_{2,1}^{\frac{n}{2}}}^2 \|\d \|_{\dot{B}_{2,1}^{1+\frac{n}{2}}}^2.
\end{align*}
The combination of  \eqref{B5}, \eqref{A2} and  \eqref{A1} implies
\begin{align}\label{A3}
&\|\delta\|_{L^\infty_T(\dot{B}_{2,1}^{\frac{n}{2}} )}+\|\delta\|_{L^1_T(\dot{B}_{2,1}^{2+\frac{n}{2}} )}\nonumber\\
&\quad\le
\varepsilon\int_0^T\|\p \u \|_{\dot{B}_{2,1}^{1+\frac{n}{2}}}dt+C
\int_0^TH_3(\p \u ,\q \u ,\v  ,\delta,\d )(\|\delta\|_{\dot{B}_{2,1}^{\frac{n}{2}}}+\|\p \u \|_{\dot{B}_{2,1}^{-1+\frac{n}{2}}})dt
\nonumber\\
&\quad\quad+\nu^{-\frac12}\left(\int_0^T\|\q \u \|_{\dot{B}_{2,1}^{-1+\frac{n}{2}}}\|\nu\nabla^2\q \u \|_{\dot{B}_{2,1}^{-1+\frac{n}{2}}}dt\right)^\frac12
\left(\int_0^T\| \d \|_{\dot{B}_{2,1}^{1+\frac{n}{2}}}^2dt\right)^\frac12
.
\end{align}

Following the derivation of \eqref{B3}, we deduce from the first equation of  \eqref{i} that
\begin{align}\label{A5}
&\|\p \u \|_{L^\infty_T(\dot{B}_{2,1}^{-1+\frac{n}{2}} )}+\|\nabla^2\p \u \|_{L^1_T(\dot{B}_{2,1}^{-1+\frac{n}{2}} )}\nonumber\\
&\quad\lesssim\int_0^T\|\nabla (\u +\v  )\|_{\dot{B}_{2,1}^{\frac{n}{2}}}\|\p \u \|_{\dot{B}_{2,1}^{-1+\frac{n}{2}}}dt
+\int_0^T\|H_1\|_{\dot{B}_{2,1}^{-1+\frac{n}{2}}}dt.
\end{align}
Now we  deal with the individual  terms of  $H_1$.
Firstly, from \eqref{daishu}, we have
\begin{align}\label{A6}
\int_0^T\|H_1^{(1)}\|_{\dot{B}_{2,1}^{-1+\frac{n}{2}}}dt&\lesssim\nu^{-1}\|(\q \u _t+\nabla a,\p \u _t,\v  _t)\|_{L^1_T(\dot{B}_{2,1}^{-1+\frac{n}{2}} )}\|\nu a \|_{L^\infty_T(\dot{B}_{2,1}^{\frac{n}{2}} )}\nonumber\\
&\lesssim\nu^{-1}(Y_n(T)+W_n(T)+{V}_n(T))X_n(T).
\end{align}
Similarly,  by \eqref{daishu} and the interpolation inequality \eqref{inter}, we have
\begin{align}\label{A10}
\|H_1^{(2)}\|_{\dot{B}_{2,1}^{-1+\frac{n}{2}}}&\lesssim(1+\|a\|_{\dot{B}_{2,1}^{\frac{n}{2}}})(\|\nabla \v  \|_{\dot{B}_{2,1}^{\frac{n}{2}}}+\|\nabla\q \u \|_{\dot{B}_{2,1}^{\frac{n}{2}}})\|\p \u \|_{\dot{B}_{2,1}^{-1+\frac{n}{2}}},\nonumber\\
\|H_1^{(3)}\|_{\dot{B}_{2,1}^{-1+\frac{n}{2}}}
&\lesssim(1+\|a\|_{\dot{B}_{2,1}^{\frac{n}{2}}})\|\q \u \|_{\dot{B}_{2,1}^{\frac{n}{2}}}\|\v  \|_{\dot{B}_{2,1}^{\frac{n}{2}}}\nonumber\\
&\lesssim(1+\|a\|_{\dot{B}_{2,1}^{\frac{n}{2}}})\|\q \u \|_{\dot{B}_{2,1}^{-1+\frac{n}{2}}}^{\frac12}\|\q \u \|_{\dot{B}_{2,1}^{1+\frac{n}{2}}}^{\frac12}\|\v  \|_{\dot{B}_{2,1}^{-1+\frac{n}{2}}}^{\frac12}\|\v  \|_{\dot{B}_{2,1}^{1+\frac{n}{2}}}^{\frac12},
\end{align}
from which we can get
\begin{align}\label{A12}
\int_0^T\|H_1^{(3)}\|_{\dot{B}_{2,1}^{-1+\frac{n}{2}}}dt
\lesssim(1+\nu^{-1} X_n(T))\nu^{-\frac12}X_n(T)^{\frac12}Y_n(T)^{\frac12}{V}_n(T).
\end{align}
By using a similar derivation of estimates \eqref{A6} and \eqref{A12}, one infer that
\begin{align}\label{A13}
\int_0^T\|H_1^{(4)}\|_{\dot{B}_{2,1}^{-1+\frac{n}{2}}}dt
&\lesssim\|a\|_{L^\infty_T(\dot{B}_{2,1}^{-1+\frac{n}{2}} )}\|\u +\v  \|_{L^\infty_T(\dot{B}_{2,1}^{-1+\frac{n}{2}} )}\|\nabla \p \u \|_{L^1_T(\dot{B}_{2,1}^{\frac{n}{2}} )}\nonumber\\
&\lesssim\nu^{-1} X_n(T)(Z_n(T)+X_n(T)+{V}_n(T))W_n(T),
\end{align}
\begin{align}\label{A15}
\int_0^T\|H_1^{(5)}\|_{\dot{B}_{2,1}^{-1+\frac{n}{2}}}dt
&\lesssim\|a\|_{L^\infty_T(\dot{B}_{2,1}^{-1+\frac{n}{2}} )}(\|\v  \|_{L^\infty_T(\dot{B}_{2,1}^{-1+\frac{n}{2}})}
\|\nabla \v  \|_{L^1_T(\dot{B}_{2,1}^{\frac{n}{2}})}) \nonumber\\
&\quad\quad+\|a\|_{L^\infty_T(\dot{B}_{2,1}^{-1+\frac{n}{2}} )}\|\q \u \|_{L^\infty_T(\dot{B}_{2,1}^{-1+\frac{n}{2}})}
\|\nabla \q \u \|_{L^1_T(\dot{B}_{2,1}^{\frac{n}{2}})})\nonumber\\
&\lesssim\nu^{-1} X_n(T)({V}_n(T)^2+\nu^{-1} X_n(T)Y_n(T)),
\end{align}
and
\begin{align}\label{A9}
\int_0^T\|
H_1^{(6)}
\|_{\dot{B}_{2,1}^{-1+\frac{n}{2}}}dt
&\lesssim\int_0^T\|\nabla \delta\odot \nabla \delta+\nabla\delta\odot \nabla \d +\nabla\d \odot \nabla \delta\|_{\dot{B}_{2,1}^{\frac{n}{2}}}dt\nonumber\\
&\lesssim\int_0^T\| \delta\|_{\dot{B}_{2,1}^{1+\frac{n}{2}}}^2+\| \d \|_{\dot{B}_{2,1}^{1+\frac{n}{2}}}\| \delta\|_{\dot{B}_{2,1}^{1+\frac{n}{2}}}dt\nonumber\\
&\lesssim\int_0^T\| \delta\|_{\dot{B}_{2,1}^{2+\frac{n}{2}}}\| \delta\|_{\dot{B}_{2,1}^{\frac{n}{2}}}+\| \d \|_{\dot{B}_{2,1}^{1+\frac{n}{2}}}\| \delta\|_{\dot{B}_{2,1}^{\frac{n}{2}}}^\frac12\| \delta\|_{\dot{B}_{2,1}^{2+\frac{n}{2}}}^\frac12dt\nonumber\\
&\lesssim\ \varepsilon\int_0^T\|\delta\|_{\dot{B}_{2,1}^{2+\frac{n}{2}}}dt+\int_0^T(\| \delta\|_{\dot{B}_{2,1}^{2+\frac{n}{2}}}+\| \d \|_{\dot{B}_{2,1}^{1+\frac{n}{2}}}^2)\| \delta\|_{\dot{B}_{2,1}^{\frac{n}{2}}}dt.
\end{align}
By the estimates \eqref{A6}-\eqref{A9}, Eq. \eqref{A5} becomes that
\begin{align}\label{A16}
&\|\p \u \|_{L^\infty_T(\dot{B}_{2,1}^{-1\!+\!\frac{n}{2}} )}\!+\!\|\nabla^2\p \u \|_{L^1_T(\dot{B}_{2,1}^{-1\!+\!\frac{n}{2}} )}\nonumber\\
&\quad\lesssim\varepsilon\| \delta\|_{L^1_T(\dot{B}_{2,1}^{2\!+\!\frac{n}{2}})}\!+\!\int_0^TH_3(\p \u ,\q \u ,\v  ,\delta,\d ) Z_n(T)dt\nonumber\\
&\quad\quad
\!+\!\nu^{-\frac12}X_n(T)^{\frac12}Y_n(T)^{\frac12}{V}_n(T)\!+\!\nu^{-1} X_n(T)(Z_n(T)\!+\!X_n(T)\!+\!{V}_n(T))W_n(T)\nonumber\\
&\quad\quad\!+\!\nu^{-1}(Y_n(T)\!+\!W_n(T)\!+\!{V}_n(T))X_n(T)\!+\!\nu^{-1} X_n(T)[{V}_n(T)^2\!+\!\nu^{-1} X_n(T)Y_n(T)].
\end{align}

Assume from now on that $X_n(T)\ll\nu$.
Summing up the estimates \eqref{A3} and \eqref{A16} and choosing $\varepsilon$ small enough  give that
\begin{align*}
&Z_n(T)+W_n(T)\nonumber\\
&\quad\lesssim\int_0^TH_3(\p \u ,\q \u ,\v  ,\delta,\d ) Z_n(T)dt\nonumber\\
&\quad\quad
+\nu^{-\frac12}X_n(T)^{\frac12}Y_n(T)^{\frac12}{V}_n(T)+\nu^{-1} X_n(T)(Z_n(T)+X_n(T)+{V}_n(T))W_n(T)\nonumber\\
&\quad\quad+\nu^{-1}(Y_n(T)+W_n(T)+{V}_n(T))X_n(T)+\nu^{-1} X_n(T)[{V}_n(T)^2+\nu^{-1} X_n(T)Y_n(T)].
\end{align*}
and so, by the  Gronwall inequality,
\begin{align}\label{A18}
Z_n(T)\!+\!W_n(T)
\lesssim&\exp\Big(\int_0^TH_3(\p \u ,\q \u ,\v  ,\delta,\d )dt\Big)\Bigg\{\nu^{-\frac12}X_n(T)^{\frac12}Y_n(T)^{\frac12}{V}_n(T)\nonumber\\
&\quad\quad\!+\!\nu^{-1} X_n(T)(Z_n(T)\!+\!X_n(T)\!+\!{V}_n(T))W_n(T)\nonumber\\
&\quad\quad\!+\!\nu^{-1}(Y_n(T)\!+\!W_n(T)\!+\!{V}_n(T))X_n(T)\!+\!\nu^{-1} X_n(T){V}_n(T)^2\Bigg\}.
\end{align}
\subsection{Estimate for the compressible part }
In fact, the compressible part  (\ref{c})  without  the term $$\diverg ( \nabla (\delta+\d )\odot \nabla (\delta+\d )-\tfrac12|\nabla (\delta+\d )|^2\I)$$ is the  same with  that  of the compressible Navier-Stokes equations  in \cite{danchin2018}.
From \eqref{daishu} and the interpolation inequality \eqref{inter}, we  have
\begin{align*}\label{}
&\|\diverg ( \nabla (\delta+\d )\odot \nabla (\delta+\d )-\tfrac12|\nabla (\delta+\d )|^2\I)\|_{\dot{B}_{2,1}^{-1+\frac{n}{2}}}\nonumber\\
&\quad\lesssim\| \nabla (\delta+\d )\odot \nabla (\delta+\d )\|_{\dot{B}_{2,1}^{\frac{n}{2}}}\nonumber\\
&\quad\lesssim\|\delta \|_{\dot{B}_{2,1}^{1+\frac{n}{2}}}^2+\|\d \|_{\dot{B}_{2,1}^{1+\frac{n}{2}}}^2\nonumber\\
&\quad\lesssim\|\delta\|_{\dot{B}_{2,1}^{\frac{n}{2}}}\|\delta\|_{\dot{B}_{2,1}^{2+\frac{n}{2}}}+\|\d \|_{\dot{B}_{2,1}^{1+\frac{n}{2}}}^2.
\end{align*}
Thus, we can get similarly from the estimate  below  \cite[Eq. (3.57)]{danchin2018} (see more details in  \cite[Step 2]{danchin2018}) that
\begin{align*}
&\|(a,\nu\nabla a,\q \u )\|_{{L}^\infty_t(\dot{B}_{2,1}^{-1+\frac{n}{2}} )}+\|(\q \u _t+\nabla a,\nu\nabla^2 \q \u ,\nu\nabla^2 a^\ell,\nabla a^h)\|_{{L}^1_t(\dot{B}_{2,1}^{-1+\frac{n}{2}} )}
\nonumber\\
&\quad\lesssim\|(a,\nu\nabla a,\q \u )(0)\|_{\dot{B}_{2,1}^{-1+\frac{n}{2}}}+\int_0^t \| (\nabla  \p \u ,\nabla  \q \u ,\nabla \v  )\|_{\dot{B}_{2,1}^{\frac{n}{2}}}\|(a,\nu\nabla a,\u )\|_{\dot{B}_{2,1}^{-1+\frac{n}{2}}}d\tau\nonumber\\
&\quad\quad+\int_0^t\|(\p \u ,\v  )\|_{\dot{B}_{2,1}^{-1+\frac{n}{2}}}\|(\nabla \p \u ,\nabla \p \v  )\|_{\dot{B}_{2,1}^{\frac{n}{2}}}
 d\tau+\int_0^t(\|\delta\|_{\dot{B}_{2,1}^{\frac{n}{2}}}\|\delta\|_{\dot{B}_{2,1}^{2+\frac{n}{2}}}+\|\d \|_{\dot{B}_{2,1}^{1+\frac{n}{2}}}^2)d\tau\nonumber\\
 &\quad\quad+\|a\|_{{L}^\infty_t(\dot{B}_{2,1}^{\frac{n}{2}} )}\int_0^t\|(\p \u ,\v  )\|_{\dot{B}_{2,1}^{-1+\frac{n}{2}}}\|\nabla \q \u \|_{\dot{B}_{2,1}^{\frac{n}{2}}}
 d\tau+\|a\|_{{L}^\infty_t(\dot{B}_{2,1}^{\frac{n}{2}} )}\|(\v  _t,\p \u _t)\|_{{L}^1_t(\dot{B}_{2,1}^{-1+\frac{n}{2}} )}
\nonumber\\&\quad\quad+\nu^{-1}(\|a^\ell\|_{{L}^\infty_t(\dot{B}_{2,1}^{-1+\frac{n}{2}} )}+\|\nu\nabla a^\ell\|_{{L}^1_t(\dot{B}_{2,1}^{\frac{n}{2}} )}
+\|\nu a^h\|_{{L}^\infty_t(\dot{B}_{2,1}^{-1+\frac{n}{2}} )}+\|a^h\|_{{L}^1_t(\dot{B}_{2,1}^{\frac{n}{2}} )}).
\end{align*}
Hence, from the Gronwall inequality, we have
\begin{align}\label{A20}
&X_n(T)\!+\!Y_n(T)\nonumber\\
&\le C \exp\left(\int_0^t \| (\nabla  \p \u ,\nabla  \q \u ,\nabla \v  )\|_{\dot{B}_{2,1}^{\frac{n}{2}}}d\tau\right)\Big( X_n(0) \!+\!({V}_n(T)\!+\!Z_n(T))({V}_n(T)\!+\!W_n(T))\nonumber\\
&\quad\!+\!\nu^{-2}X_n(T)Y_n(T)({V}_n(T)\!+\!Z_n(T))\!+\!X_n(T)\!+\!\nu^{-1}({V}_n(T)\!+\!Y_n(T)\!+\!W_n(T))X_n(T)         \Big).
\end{align}

\subsection{Completion for the proof of  Theorem \ref{dingli}}

We claim that if $\nu$ is large enough then one may find a large ${\gamma}$ and a small $\eta$
so that for all $T<T^*,$  the following bounds are  valid:
\begin{equation}\label{f3}
 X_n(T)+Y_n(T) \leq {\gamma}\quad\hbox{and}\quad
 Z_n(T)+W_n(T) \leq \eta.
\end{equation}

Assuming that
\begin{align}\label{C1}
\nu^{-1}{\gamma}\ll1,
\end{align}

one deduce from \eqref{A18} and \eqref{A20} that
\begin{align}\label{A21}
X_n(T)+Y_n(T)\le&C \exp\big(M+\nu^{-1}{\gamma}+\eta\big)\Big( X_n(0) +(M+\eta)^2\nonumber\\
&\quad+\nu^{-2}{\gamma}(M+\eta)X_n(T)+\nu^{-1}(M+{\gamma}+\eta)X_n(T)
\Big),
\end{align}
\begin{align*}
Z_n(T)+W_n(T)
\le&C\exp\Big(\int_0^TH_3(\p \u ,\q \u ,\v  ,\delta,\d )dt\Big)\Bigg\{\nu^{-\frac12}X_n(T)^{\frac12}Y_n(T)^{\frac12}{V}_n(T)\nonumber\\
&\quad\quad+\nu^{-1} X_n(T)(Z_n(T)+X_n(T)+{V}_n(T))W_n(T)\nonumber\\
&\quad\quad+\nu^{-1}(Y_n(T)+W_n(T)+{V}_n(T))X_n(T)+\nu^{-1} X_n(T){V}_n(T)^2\Bigg\}.
\end{align*}
A simple computation gives
\begin{align*}
&\int_0^TH_3(\p \u ,\q \u ,\v  ,\delta,\d )dt\nonumber\\
&\quad\lesssim W_n(T)+\nu^{-1}Y_n(T)+{V}_n(T)+Z_n(T)W_n(T)+{V}_n(T)^2+{V}_n(T)W_n(T)+{V}_n(T)^4\nonumber\\
&\quad\lesssim \eta+\nu^{-1} {\gamma}+M\eta+M^2+M^4,
\end{align*}
and thus,
\begin{align}\label{A28}
Z_n(T)+W_n(T)
\le&C\exp\big(\eta+\nu^{-1} {\gamma}+M\eta+M^2+M^4\big)
\Big\{\nu^{-\frac12} {\gamma}M\nonumber\\
&+\nu^{-1} {\gamma}(\eta+{\gamma}+M)W_n(T)+\nu^{-1}({\gamma}+\eta+M){\gamma}+\nu^{-1} {\gamma}M^2\Big\}\nonumber\\
\le&C{\gamma}\exp\big(\eta+\nu^{-1} {\gamma}+M\eta+M^2+M^4\big)
\Big\{\nu^{-\frac12} M\nonumber\\
&+\nu^{-1} (\eta+{\gamma}+M)W_n(T)+\nu^{-1}({\gamma}+\eta+M)+\nu^{-1} M^2\Big\}.
\end{align}
Assumption in addition that
\begin{align}\label{C2}
\nu^{-1}{\gamma}\ll M,\ \eta \le \min\{M,1\},
\end{align}
one can deduce from \eqref{A21} and \eqref{A28} that
\begin{eqnarray}\label{mhd1}
\left\{\begin{aligned}
X_n(T)\!+\!Y_n(T)\!\le&C \exp\big(CM\big)\Big( X_n(0) \!+\!M^2\!+\!1\!+\!\nu^{-1}(M\!+\!{\gamma})X_n(T)
\Big),\\
Z_n(T)\!+\!W_n(T)
\!\le& C{\gamma}\exp\big(CM\!+\!CM^4\big)
\Big\{\nu^{-\frac12} M
\!+\!\nu^{-1} ({\gamma}\!+\!M)W_n(T)\!+\!\nu^{-1}({\gamma}\!+\!M^2)\Big\}.
\end{aligned}\right.\hspace{5mm}
\end{eqnarray}
Making further assumption that
\begin{align}\label{C3}
\nu^{-1}{\gamma}(M+{\gamma}+1)\exp\big(CM+CM^4\big)\ll 1,
\end{align}
then \eqref{mhd1} becomes that
\begin{eqnarray}
\left\{\begin{aligned}
X_n(T)+Y_n(T)\le&C \exp\big(CM\big)\Big( X_n(0) +M^2+1
\Big),\nonumber\\
Z_n(T)+W_n(T)
\le&C{\gamma}\exp\big(CM+CM^4\big)
\Big\{\nu^{-\frac12}M+
\nu^{-1}( 1
+{\gamma}+M^2)\Big\}.
\end{aligned}\right.\nonumber
\end{eqnarray}
We first define
\begin{align}\label{C5}
&{\gamma}= C\exp\big(CM+CM^4\big)\Big( X_n(0) +M^2+1\Big),
\end{align}
an then let
\begin{align}\label{C6}
&\eta= C\exp\big(CM+CM^4\big)\Big( X_n(0) +M^2
\Big)
\Big\{\nu^{-\frac12}M+\nu^{-1} \Big( X_n(0) +M^2+1
\Big)\Big\}.
\end{align}
For a suitable large constant $C$, taking
\begin{align}\label{C8}
C\exp\big(C(M+M^4)\big)\Big( X_d(0) +M^2+1
\Big)\le\sqrt{\nu},
\end{align}
then the conditions in \eqref{C1}, \eqref{C2}, \eqref{C3} are all satisfied.
Thus,  if $\nu$ and the compressible part of the data
fulfill  \eqref{C8}  then defining $\gamma$ and $\eta$ according to \eqref{C5} and \eqref{C6}
ensures that \eqref{f3} is fulfilled for all $T<T^*.$

Finally, using the fact that
$
v=\u +\v=\p \u +\q \u  +\v $, $d=\delta+\d$  and
 the interpolation inequality \eqref{inter}, we have
\begin{align*}
\int_0^{T^\ast}(\|\nabla v\|_{\dot{B}_{2,1}^{\frac{n}{2}}}+\|d\|^2_{\dot{B}_{2,1}^{1+\frac{n}{2}}})dt
\le& \int_0^{T^\ast}(\|\nabla \p \u \|_{\dot{B}_{2,1}^{\frac{n}{2}}}+\|\nabla\q \u \|_{\dot{B}_{2,1}^{\frac{n}{2}}} +\|\nabla \v \|_{\dot{B}_{2,1}^{\frac{n}{2}}})dt\nonumber\\
&+\int_0^{T^\ast}(\|\delta\|_{\dot{B}_{2,1}^{1+\frac{n}{2}}}^2+\|\d\|^2_{\dot{B}_{2,1}^{1+\frac{n}{2}}})dt\nonumber\\
\le& \int_0^{T^\ast}(\|\nabla \p \u \|_{\dot{B}_{2,1}^{\frac{n}{2}}}+\|\nabla\q \u \|_{\dot{B}_{2,1}^{\frac{n}{2}}} +\|\nabla \v \|_{\dot{B}_{2,1}^{\frac{n}{2}}})dt\nonumber\\
&+\int_0^{T^\ast}(\|\delta\|_{\dot{B}_{2,1}^{\frac{n}{2}}}\|\delta\|_{\dot{B}_{2,1}^{2+\frac{n}{2}}}+\|\d\|_{\dot{B}_{2,1}^{\frac{n}{2}}}\|\d\|_{\dot{B}_{2,1}^{2+\frac{n}{2}}})dt
\nonumber\\
\le& C(\eta+\gamma\nu^{-1}+M+M^2+\eta^2)<\infty.
\end{align*}
Thus,  combining with the
continuation criterion recalled in \eqref{baopobiaozhun}, one can conclude that $T^*=+\infty$.

Consequently, we have completed the proof of  Theorem \ref{dingli}.

\appendix

 \section{Proof of the blow-up criterion \eqref{baopobiaozhun} in Theorem \ref{dingli2}}

We prove \eqref{baopobiaozhun} by contradiction. Let $0 < T^\ast < \infty$ be the maximum time for the existence
of strong solution $(a, v, d)$ to the system (\ref{m}). Assume that \eqref{baopobiaozhun} is not true. Then there exists
a positive constant $C_0$ such that

\begin{align}\label{baopobiaozhun1}
\int_0^{T^\ast}(\|\nabla v\|_{\dot{B}_{2,1}^{\frac{n}{2}}}+\|d\|^2_{\dot{B}_{2,1}^{1+\frac{n}{2}}})dt\le C_0.
\end{align}
Under the assumption \eqref{baopobiaozhun1},  it thus remains  to show the existence of a  bound $C$ depending only on $a_0,v_0, d_0, T^\ast$
and $C_0$ such that
\begin{align}\label{}
\sup_{0<t\le T^\ast}(\|a\|_{\widetilde{L}_t^{\infty}(\dot{B}_{2,1}^{\frac n2})}+\|v\|_{\widetilde{L}_t^{\infty}(\dot{B}_{2,1}^{-1+\frac n2})}
+\|d\|_{\widetilde{L}_t^{\infty}(\dot{B}_{2,1}^{\frac n2})})\le C.
\end{align}
To do so, we get from the first equation in \eqref{m} that
\begin{align*}\label{}
\partial_t a+ v\cdot\nabla a+(1+a)\diverg v=0.
\end{align*}
Applying $\dot{\Delta}_j$ to the previous  equation gives 
\begin{align*}\label{}
\partial_t \dot{\Delta}_ja+v\cdot\nabla \dot{\Delta}_j a+\dot{\Delta}_j((1+a)\diverg v)=[v\cdot\nabla, \dot{\Delta}_j ]a.
\end{align*}
Taking $L^2$ inner product with $ \dot{\Delta}_j a $ and then integrating in time on $[0,T^*]$, we have
\begin{align*}
\|\dot{\Delta}_ja\|_{L^2}\le \|\dot{\Delta}_ja_0\|_{L^2}\!+\!\int_0^{T^*}\Big( \|\diverg v\|_{L^\infty}\|\dot{\Delta}_ja\|_{L^2}
\!+\!\|\dot{\Delta}_j((1+a)\diverg v)\|_{L^2}\!+\!\|[\dot{\Delta}_j,v\cdot\nabla  ]a\|_{L^2}\Big)
dt.
\end{align*}
Multiplying by $2^{\frac n2j}$ to the previous equation  and then summing up the resultant equation with respect to  $j\in\Z$, we can get from \eqref{daishu} and \eqref{jiaohuanzi} that
\begin{align*}\label{}
\|a\|_{\dot{B}_{2,1}^{\frac{n}{2}}}\le\|a_0\|_{\dot{B}_{2,1}^{\frac{n}{2}}}+\int_0^{T^*}\|\nabla v\|_{\dot{B}_{2,1}^{\frac{n}{2}}}dt
+\int_0^{T^*}\|\nabla v\|_{\dot{B}_{2,1}^{\frac{n}{2}}}\|a\|_{\dot{B}_{2,1}^{\frac{n}{2}}}dt.
\end{align*}
This, together with the Gronwall inequality, implies that
\begin{align}\label{m1}
\|a\|_{\widetilde{L}^\infty_{T^\ast}({\dot{B}_{2,1}^{\frac n2}})}\le&(\|a_0\|_{\dot{B}_{2,1}^{\frac{n}{2}}}+\int_0^{T^*}\|\nabla v\|_{\dot{B}_{2,1}^{\frac{n}{2}}}dt)\exp{(C\|\nabla v\|_{L_{T^\ast}^1({\dot{B}_{2,1}^{\frac{n}{2}}})})}\nonumber\\
\le& \exp(CC_0)(C_0+\|a_0\|_{\dot{B}_{2,1}^{\frac{n}{2}}})<\infty.
\end{align}
Similarly, we can get from the second equation in \eqref{m} that
\begin{align*}
&\|d\|_{\widetilde{L}^\infty_{T^\ast}({\dot{B}_{2,1}^{\frac n2}})}+\|d\|_{L^1_{T^\ast}({\dot{B}_{2,1}^{2+\frac n2}})}\nonumber\\
&\quad\lesssim\|d_0\|_{\dot{B}_{2,1}^{\frac n2}}+\int_0^{T^\ast}(\|\nabla \u \|_{\dot{B}_{2,1}^{\frac{n}{2}}}+\|d\|^2_{\dot{B}_{2,1}^{1+\frac{n}{2}}} )\|d\|_{\dot{B}_{2,1}^{\frac{n}{2}}}dt+\int_0^{T^\ast}\|d\|^2_{\dot{B}_{2,1}^{1+\frac{n}{2}}} dt
\end{align*}
and hence, by the Gronwall inequality,
\begin{align}\label{m3}
&\|d\|_{\widetilde{L}^\infty_{T^\ast}({\dot{B}_{2,1}^{\frac n2}})}+\|d\|_{L^1_{T^\ast}({\dot{B}_{2,1}^{2+\frac n2}})}\nonumber\\
&\quad\lesssim(\|d_0\|_{\dot{B}_{2,1}^{\frac n2}}+\int_0^{T^\ast}\|d\|^2_{\dot{B}_{2,1}^{1+\frac{n}{2}}} dt)\exp(\int_0^{T^\ast}(\|\nabla v\|_{\dot{B}_{2,1}^{\frac{n}{2}}}+\|d\|^2_{\dot{B}_{2,1}^{1+\frac{n}{2}}} )dt)\nonumber\\
&\quad\lesssim(\|d_0\|_{\dot{B}_{2,1}^{\frac n2}}+C_0)\exp(C_0)\nonumber\\
&\quad<\infty.
\end{align}
Finally, to show the estimate of $v$, we   deduce from the second equation in \eqref{m} that
\begin{align}
v_t&+v\cdot\nabla v-\frac{\mu}{1+a} \Delta v-\frac{\lambda+\mu}{1+a}\nabla\diverg v\nonumber\\
&=-\frac{1}{1+a}\big(\nabla P(1+a)+\diverg ( \nabla d\odot \nabla d-\tfrac12|\nabla d|^2\I)\big).
\label{fangcheng}
\end{align}
By the estimate
$\|a\|_{\widetilde{L}^\infty_{T^\ast}({\dot{B}_{2,1}^{\frac n2}})}<\infty$ described in \eqref{m1}, we can find an integer  $m\in\Z$  large enough such that
\begin{align*}
&\inf_{(t,x)\in[0,T^*)\times \R^n}(1+\dot{S}_m a(t,x))\ge \frac12,\quad
\|a-\dot{S}_m a\|_{{L}^\infty_{T^\ast}({\dot{B}_{2,1}^{\frac n2}})}\le {c\underline{b}}{\overline{b}^{-1}},
\end{align*}
with
$
\overline{b}= \mu+|\lambda+\mu|$ and $\underline{b}= \min(\mu,\lambda+2\mu).
$

Thus, applying
 \cite[Proposition 6]{danchin2007cpde} to  \eqref{fangcheng} gives
\begin{align}\label{m5}
&\|v\|_{\widetilde{L}^\infty_{T^\ast}({\dot{B}_{2,1}^{-1+\frac n2}})}+\|v\|_{L^1_{T^\ast}({\dot{B}_{2,1}^{1+\frac n2}})}\nonumber\\
&\quad\le  C\exp{\Big(\int_0^{T^\ast}C(\|\nabla v\|_{\dot{B}_{2,1}^{\frac{n}{2}}}+2^{2m}\bar{b}^2b^{-1}\|a\|_{\dot{B}_{2,1}^{\frac n2}}^2)dt
\Big)}\nonumber\\
&\quad\quad\times\Big(\|v_0\|_{\dot{B}_{2,1}^{-1+\frac n2}}+\int_0^{T^\ast}\Big\|\frac{a}{1+a}(\nabla a+\diverg ( \nabla d\odot \nabla d-\tfrac12|\nabla d|^2\I))\Big\|_{\dot{B}_{2,1}^{-1+\frac n2}}dt\Big)\nonumber\\
&\quad\le  C\exp{\Big(\int_0^{T^\ast}C(\|\nabla v\|_{\dot{B}_{2,1}^{\frac{n}{2}}}+2^{2m}\overline{{b}}^2\underline{b}^{-1}\|a\|_{\dot{B}_{2,1}^{\frac n2}}^2)dt
\Big)}\nonumber\\
&\quad\quad\times\Big(\|v_0\|_{\dot{B}_{2,1}^{-1+\frac n2}}+\int_0^{T^\ast}(1+\|a\|_{\dot{B}_{2,1}^{\frac n2}})(\|a\|_{\dot{B}_{2,1}^{\frac n2}}+\|d\|^2_{\dot{B}_{2,1}^{1+\frac{n}{2}}}) dt\Big)\nonumber\\
&\quad\le  C\exp{\Big(CC_0+C2^{2m}\overline{{b}}^2\underline{b}^{-1}T^*\exp(2CC_0)(C_0+\|a_0\|_{\dot{B}_{2,1}^{\frac{n}{2}}})^2
\Big)}\nonumber\\
&\quad\quad\times\left(\|v_0\|_{\dot{B}_{2,1}^{-1+\frac n2}}+\Big(1+\exp(CC_0)(C_0+\|a_0\|_{\dot{B}_{2,1}^{\frac{n}{2}}})\Big)\Big(\exp(CC_0)(C_0+
\|a_0\|_{\dot{B}_{2,1}^{\frac{n}{2}}})+C_0\Big)\right)\nonumber\\
&\quad<\infty.
\end{align}
Therefore, combining  \eqref{m1}, \eqref{m3} and \eqref{m5}, we can extend the solution $(a, v, d)$ beyond $T^\ast$, which leads to a contradiction. The proof is complete.

\bigskip

\bigskip
\noindent \textbf{Acknowledgement.} This work is supported by NSFC under grant numbers 11601533 and 11571240.


\newpage


\noindent \textbf{References}
\begin{thebibliography}{10}

\bibitem{bcd}
H.~Bahouri, J.~Y. Chemin, R.~Danchin.
\newblock {\it {F}ourier {A}nalysis and {N}onlinear {P}artial {D}ifferential
  {E}quations}.
\newblock Grundlehren Math. Wiss., vol. {\textbf{343}}, Springer-Verlag,
  Berlin, Heidelberg, 2011.
\newblock $\,$
\bibitem{bai1}
Q. Bie, H. Cui, Q. Wang, Z. Yao. Global existence and incompressible limit in critical
spaces for compressible
ow of liquid crystals.  {\it Z. Angew. Math. Phys.}, {\bf 68} (2017), 113.

\bibitem{bai2}
Q. Bie, H. Cui, Q. Wang, Z. Yao. Incompressible limit for the compressible flow of liquid crystals in $L^p$ critical Besov spaces.  {\it Discrete Contin. Dyn. Syst.}, {\bf 38} (2018), 2879--2910.



\bibitem{danchin2007cpde}
R.~Danchin.
\newblock Well-posedness in critical spaces for barotropic viscous fluids with
  truly not constant density.
\newblock {\it Comm. Partial Differential Equations}, {\bf 32} (2007), 1373--1397.
\newblock $\,$

\bibitem{danchin2000} R. Danchin. Global existence in critical spaces for compressible
Navier-Stokes equations. {\it Invent. Math.}, {\bf 141} (2000),  579--614.



\bibitem{danchin2018}
R.~Danchin, P. Mucha.
\newblock
Compressible Navier-Stokes system: large solutions and incompressible limit.
\newblock {\it Adv. Math.}, {\bf 320} (2017), 904--925.
\newblock $\,$



\bibitem{ding1} S. Ding, J. Huang, H. Wen, R. Zi.
Incompressible limit of the compressible nematic liquid crystal flow. {\it J. Funct. Anal.}, {\bf 264} (2013),  1711--1756.



\bibitem{hao}
Y. Hao, X. Liu.
\newblock Incompressible limit of a compressible liquid crystals system.
\newblock {\it Acta Math. Sci. Ser. B Engl. Ed.}, {\bf 33} (2013), 781--796.

\bibitem{Hoff}
D. Hoff, K. Zumbrun.  Multidimensional diffusion waves for the
Navier-Stokes equations of compressible flow. { \it Indiana Univ. Math. J.}, {\bf{44}} (1995),
604--676.

\bibitem{huxianpeng}
X.~Hu,  H.~Wu.
\newblock Global solution to the three-dimensional compressible flow of liquid crystals.
\newblock {\it  SIAM J. Math. Anal.}, {\bf 45} (2013), 2678--2699.
\newblock $\,$




\bibitem{huang2016} T. Huang, F. Lin, C. Liu,  C. Wang.
Finite time singularity of the nematic liquid crystal flow in dimension three. {\it Arch. Ration. Mech. Anal.}, {\bf 221} (2016), 1223--1254.

\bibitem{huang2012}T. Huang, C.  Wang, H.  Wen.
Strong solutions of the compressible nematic liquid crystal flow. {\it J. Differential Equations}, {\bf 252} (2012), 2222--2265.

\bibitem{huang2012+1} T. Huang, C. Wang, H. Wen.
Blow up criterion for compressible nematic liquid crystal flows in dimension three. {\it Arch. Ration. Mech. Anal.}, {\bf 204} (2012), 285--311.


\bibitem{huangxiangdi2012}
X. Huang, J. Li, Z.  Xin. Global well-posedness of classical
solutions with large oscillations and vacuum to the
three-dimensional isentropic compressible Navier-Stokes equaitons.
{\it Commun. Pure Appl. Math.}, {\bf{65}} (2012), 549--585.





\bibitem{jjw} F. Jiang. S. Jiang, D. Wang.
On multi-dimensional compressible flows of nematic liquid crystals with large initial energy in a bounded domain. {\it J. Funct. Anal.}, {\bf 265} (2013),  3369--3397.

\bibitem{jiang2014}
F. Jiang, S. Jiang,  D. Wang.
\newblock Global weak solutions to the equations of compressible flow of
  nematic liquid crystals in two dimensions.
\newblock {\it Arch. Ration. Mech. Anal.},  {\bf 214} (2014), 403--451.





\bibitem{linfanghua2016} F. Lin,  C. Wang.
Global existence of weak solutions of the nematic liquid crystal flow in dimension three. {\it Comm. Pure Appl. Math.}, {\bf 69} (2016),  1532--1571.

\bibitem{linfanghua2015} J. Lin, B. Lai, C. Wang.
Global finite energy weak solutions to the compressible nematic liquid crystal flow in dimension three. {\it SIAM J. Math. Anal.}, {\bf 47} (2015), 2952--2983.

\bibitem{lions}
P. L. Lions, N. Masmoudi.
\newblock Incompressible limit for a viscous compressible fluid.
\newblock {\it J. Math.  Pures  Appl.}
  {\bf 77} (1998), 585--627.


\bibitem{liuqiao}
Q.~Liu,  T.~Zhang, J.~Zhao.
\newblock Global solutions to the 3D incompressible nematic liquid crystal system.
\newblock {\it J. Differential Equations}, {\bf 258} (2015), 1519--1547.
\newblock $\,$

\bibitem{liu2013} X. Liu,  J. Qing.
Globally weak solutions to the flow of compressible liquid crystals system. {\it Discrete Contin. Dyn. Syst.}, {\bf 33} (2013), 757--788.



\bibitem{ouraobing}
Y. Ou.
\newblock Low mach number limit of viscous polytropic fluid flows.
\newblock {\it J. Differential Equations}, {\bf 251} (2011), 2037--2065.

\bibitem{WangA2012} D. Wang, C. Yu.
Global weak solution and large-time behavior for the compressible flow of liquid crystals. {\it Arch. Ration. Mech. Anal.}, {\bf 204} (2012), 881--915.

\bibitem{wang2014}
D. Wang, C. Yu.
\newblock Incompressible limit for the compressible flow of liquid crystals.
\newblock {\it J. Math. Fluid Mech.}, {\bf 16} (2014), 771--786.


\bibitem{zhaixiaoping}
X.~Zhai,  Z.~Yin.
\newblock
On some large global solutions to the incompressible inhomogeneous nematic liquid
crystal flows,
\newblock {\it Accepted for  publication}.
\newblock $\,$

















\end{thebibliography}
\end{document}